\documentclass[11pt]{article}

\usepackage{latexsym}
\usepackage[dvips]{graphicx}
\usepackage{float}
\usepackage{amsmath}
\usepackage{mathdots}   
\usepackage{amssymb}
\usepackage{epic}
\usepackage{color}
\usepackage{setspace}
\usepackage{lscape}
\usepackage{arydshln}
\usepackage{fancybox}   
\usepackage{mathtools}
\usepackage{amsthm}

\begin{document}
\bibliographystyle{plain}
\floatplacement{table}{H}
\newtheorem{definition}{Definition}[section]
\newtheorem{lemma}{Lemma}[section]
\newtheorem{theorem}{Theorem}[section]
\newtheorem{corollary}{Corollary}[section]
\newtheorem{proposition}{Proposition}[section]

\newcommand{\sni}{\sum_{i=1}^{n}}
\newcommand{\snj}{\sum_{j=1}^{n}}
\newcommand{\smj}{\sum_{j=1}^{m}}
\newcommand{\sumjm}{\sum_{j=1}^{m}}
\newcommand{\bdis}{\begin{displaymath}}
\newcommand{\edis}{\end{displaymath}}
\newcommand{\beq}{\begin{equation}}
\newcommand{\eeq}{\end{equation}}
\newcommand{\beqn}{\begin{eqnarray}}
\newcommand{\eeqn}{\end{eqnarray}}
\newcommand{\defeq}{\stackrel{\triangle}{=}}
\newcommand{\simleq}{\stackrel{<}{\sim}}
\newcommand{\sep}{\;\;\;\;\;\; ; \;\;\;\;\;\;}
\newcommand{\real}{\mbox{$ I \hskip -4.0pt R $}}
\newcommand{\complex}{\mbox{$ I \hskip -6.8pt C $}}
\newcommand{\integ}{\mbox{$ Z $}}
\newcommand{\realn}{\real ^{n}}
\newcommand{\sqrn}{\sqrt{n}}
\newcommand{\sqrtwo}{\sqrt{2}}
\newcommand{\prf}{{\bf Proof. }}

\newcommand{\onehlf}{\frac{1}{2}}
\newcommand{\thrhlf}{\frac{3}{2}}
\newcommand{\fivhlf}{\frac{5}{2}}
\newcommand{\onethd}{\frac{1}{3}}
\newcommand{\lb}{\left ( }
\newcommand{\lcb}{\left \{ }
\newcommand{\lsb}{\left [ }
\newcommand{\labs}{\left | }
\newcommand{\rb}{\right ) }
\newcommand{\rcb}{\right \} }
\newcommand{\rsb}{\right ] }
\newcommand{\rabs}{\right | }
\newcommand{\lnm}{\left \| }
\newcommand{\rnm}{\right \| }
\newcommand{\lambdab}{\bar{\lambda}}
%
%
\newcommand{\xj}{x_{j}}
\newcommand{\xjb}{\bar{x}_{j}}
\newcommand{\xro}{x_{\resh}}
\newcommand{\xrob}{\bar{x}_{\resh}}
\newcommand{\xsig}{x_{\sigma}}
\newcommand{\xsigb}{\bar{x}_{\sigma}}
\newcommand{\xnmjb}{\bar{x}_{n-j+1}}
\newcommand{\xnmj}{x_{n-j+1}}
\newcommand{\aroj}{a_{\resh j}}
\newcommand{\arojb}{\bar{a}_{\resh j}}
\newcommand{\aroro}{a_{\resh \resh}}
\newcommand{\amuro}{a_{\mu \resh}}
\newcommand{\amumu}{a_{\mu \mu}}
\newcommand{\aii}{a_{ii}}
\newcommand{\aik}{a_{ik}}
\newcommand{\akj}{a_{kj}}
\newcommand{\atwoii}{a^{(2)}_{ii}}
\newcommand{\atwoij}{a^{(2)}_{ij}}
\newcommand{\ajj}{a_{jj}}
\newcommand{\aiib}{\bar{a}_{ii}}
\newcommand{\ajjb}{\bar{a}_{jj}}
\newcommand{\bii}{a_{jj}}
\newcommand{\biib}{\bar{a}_{jj}}
\newcommand{\aij}{a_{i,n-i+1}}
\newcommand{\akl}{a_{j,n-j+1}}
\newcommand{\aijb}{\bar{a}_{i,n-i+1}}
\newcommand{\aklb}{\bar{a}_{j,n-j+1}}
\newcommand{\bij}{a_{n-j+1,j}}
\newcommand{\arorob}{\bar{a}_{\resh \resh}}
\newcommand{\arosig}{a_{\resh \sigma}}
\newcommand{\arosigb}{\bar{a}_{\resh \sigma}}
\newcommand{\sumjrosig}{\sum_{\stackrel{j=1}{j\neq\resh,\sigma}}^{n}}
\newcommand{\summuro}{\sum_{\stackrel{j=1}{j\neq\mu,\resh}}^{n}}
\newcommand{\sumjnoti}{\sum_{\stackrel{j=1}{j\neq i}}^{n}}
\newcommand{\sumlnoti}{\sum_{\stackrel{\ell=1}{\ell \neq i}}^{n}}
\newcommand{\sumknoti}{\sum_{\stackrel{k=1}{k\neq i}}^{n}}
\newcommand{\sumknotij}{\sum_{\stackrel{k=1}{k\neq i,j}}^{n}}
\newcommand{\sumk}{\sum_{k=1}^{n}}
\newcommand{\snl}{\sum_{\ell=1}^{n}}
\newcommand{\sumji}{\sum_{\stackrel{j=1}{j\neq i, n-i+1}}^{n}}
\newcommand{\sumki}{\sum_{\stackrel{k=1}{k\neq i, n-i+1}}^{n}}
\newcommand{\sumkj}{\sum_{\stackrel{k=1}{k\neq j, n-j+1}}^{n}}
\newcommand{\sumjro}{\sum_{\stackrel{j=1}{j\neq\resh}}^{n}}
\newcommand{\rrosig}{R''_{\resh \sigma}}
\newcommand{\rro}{R'_{\resh}}
\newcommand{\gamror}{\Gamma_{\resh}^{R}(A)}
\newcommand{\gamir}{\Gamma_{i}^{R}(A)}
\newcommand{\gamctrr}{\Gamma_{\frac{n+1}{2}}^{R}(A)}
\newcommand{\gamctrc}{\Gamma_{\frac{n+1}{2}}^{C}(A)}
\newcommand{\gamroc}{\Gamma_{\resh}^{C}(A)}
\newcommand{\gamjc}{\Gamma_{j}^{C}(A)}
\newcommand{\lamror}{\Lambda_{\resh}^{R}(A)}
\newcommand{\lamir}{\Lambda_{i}^{R}(A)}
\newcommand{\lamirepsilon}{\Lambda_{i}^{R}(A_{\epsilon})}
\newcommand{\lamnir}{\Lambda_{n-i+1}^{R}(A)}
\newcommand{\lamjr}{\Lambda_{j}^{R}(A)}
\newcommand{\varphiij}{\Phi_{ij}^{R}(A)}
\newcommand{\delir}{\Delta_{i}^{R}(A)}
\newcommand{\vir}{V_{i}^{R}(A)}
\newcommand{\pamir}{\Pi_{i}^{R}(A)}
\newcommand{\xir}{\Xi_{i}^{R}(A)}
\newcommand{\lamjc}{\Lambda_{j}^{C}(A)}
\newcommand{\vjc}{V_{j}^{C}(A)}
\newcommand{\pamjc}{\Pi_{j}^{C}(A)}
\newcommand{\xjc}{\Xi_{j}^{C}(A)}
\newcommand{\lamroc}{\Lambda_{\resh}^{C}(A)}
\newcommand{\lamsigr}{\Lambda_{\sigma}^{R}(A)}
\newcommand{\lamsigc}{\Lambda_{\sigma}^{C}(A)}
\newcommand{\psii}{\Psi_{i}^{R}(A)}
\newcommand{\psiq}{\Psi_{q}}
\newcommand{\psiiepsilon}{\Psi_{i}^{R}(A_{\epsilon})}
\newcommand{\psiqepsilon}{\Psi_{q}(A_{\epsilon})}
\newcommand{\psiqc}{\Psi_{q}^{c}}
\newcommand{\psiqcepsilon}{\Psi_{q}^{c}(A_{\epsilon})}

\newcommand{\xmu}{x_{\mu}}
\newcommand{\xmub}{\bar{x}_{\mu}}
\newcommand{\xnu}{x_{\nu}}
\newcommand{\xnub}{\bar{x}_{\nu}}
\newcommand{\amuj}{a_{\mu j}}
\newcommand{\amujb}{\bar{a}_{\mu j}}
\newcommand{\amumub}{\bar{a}_{\mu \mu}}
\newcommand{\amunu}{a_{\mu \nu}}
\newcommand{\amunub}{\bar{a}_{\mu \nu}}
\newcommand{\sumjmunu}{\sum_{\stackrel{j=1}{j\neq\mu,\nu}}}
\newcommand{\rmunu}{R''_{\mu \nu}}
\newcommand{\rmu}{R'_{\mu}}

\newcommand{\Mzero}{M_{0}}
\newcommand{\Mone}{M_{1}}
\newcommand{\Mtwo}{M_{2}}
\newcommand{\Mth}{M_{3}}
\newcommand{\Mfr}{M_{4}}
\newcommand{\Mfv}{M_{5}}
\newcommand{\Msx}{M_{6}}
\newcommand{\Mnmo}{M_{n-1}}
\newcommand{\Mnpo}{M_{n+1}}
\newcommand{\Mnmt}{M_{n-2}}
\newcommand{\Mn}{M_{n}}
\newcommand{\Mj}{M_{j}}
\newcommand{\Mk}{M_{k}}

\newcommand{\Nzero}{N_{0}}
\newcommand{\None}{N_{1}}
\newcommand{\Ntwo}{N_{2}}
\newcommand{\Nth}{N_{3}}
\newcommand{\Nfr}{N_{4}}
\newcommand{\Nfv}{N_{5}}
\newcommand{\Nsx}{N_{6}}
\newcommand{\Nnmo}{N_{n-1}}
\newcommand{\Nnpo}{N_{n+1}}
\newcommand{\Nnmt}{N_{n-2}}
\newcommand{\Nn}{N_{n}}
\newcommand{\Nj}{N_{j}}

\newcommand{\Dzero}{D_{0}}
\newcommand{\Done}{D_{1}}
\newcommand{\Dtwo}{D_{2}}
\newcommand{\Dth}{D_{3}}
\newcommand{\Dfr}{D_{4}}
\newcommand{\Dfv}{D_{5}}
\newcommand{\Dsx}{D_{6}}
\newcommand{\Dnmo}{D_{n-1}}
\newcommand{\Dnpo}{D_{n+1}}
\newcommand{\Dnmt}{D_{n-2}}
\newcommand{\Dn}{D_{n}}
\newcommand{\Dj}{D_{j}}

\newcommand{\Azero}{A_{0}}
\newcommand{\Aone}{A_{1}}
\newcommand{\Atwo}{A_{2}}
\newcommand{\Ath}{A_{3}}
\newcommand{\Afr}{A_{4}}
\newcommand{\Afv}{A_{5}}
\newcommand{\Asx}{A_{6}}
\newcommand{\Anmo}{A_{n-1}}
\newcommand{\Anpo}{A_{n+1}}
\newcommand{\Anmt}{A_{n-2}}
\newcommand{\An}{A_{n}}
\newcommand{\Aj}{A_{j}}

\newcommand{\azero}{a_{0}}
\newcommand{\aone}{a_{1}}
\newcommand{\atwo}{a_{2}}
\newcommand{\ath}{a_{3}}
\newcommand{\afr}{a_{4}}
\newcommand{\afv}{a_{5}}
\newcommand{\asx}{a_{6}}
\newcommand{\anmo}{a_{n-1}}
\newcommand{\anmt}{a_{n-2}}
\newcommand{\an}{a_{n}}
\newcommand{\aj}{a_{j}}

\newcommand{\Bzero}{B_{0}}
\newcommand{\Bone}{B_{1}}
\newcommand{\Btwo}{B_{2}}
\newcommand{\Bth}{B_{3}}
\newcommand{\Bfr}{B_{4}}
\newcommand{\Bfv}{B_{5}}
\newcommand{\Bsx}{B_{6}}
\newcommand{\Bnmo}{B_{n-1}}
\newcommand{\Bndtwomo}{B_{n/2-1}}
\newcommand{\Bnmt}{B_{n-2}}
\newcommand{\Bn}{B_{n}}
\newcommand{\Bj}{B_{j}}

\newcommand{\Anmk}{A_{n-k}}
\newcommand{\Anmi}{A_{n-i}}
\newcommand{\Bnmk}{A_{n-k}}
\newcommand{\anmk}{a_{n-k}}
\newcommand{\zj}{z^{j}}          
\newcommand{\zk}{z^{k}}          

\newcommand{\cii}{c_{ii}}
\newcommand{\cik}{c_{ik}}
\newcommand{\ckj}{c_{kj}}
\newcommand{\ctwoii}{c^{(2)}_{ii}}
\newcommand{\ctwoij}{c^{(2)}_{ij}}
\newcommand{\cjj}{c_{jj}}

\newcommand{\bik}{b_{ik}}
\newcommand{\bkj}{b_{kj}}
\newcommand{\btwoii}{b^{(2)}_{ii}}
\newcommand{\btwoij}{b^{(2)}_{ij}}
\newcommand{\bjj}{b_{jj}}

\newcommand{\abii}{(AB)_{ii}}
\newcommand{\abil}{(AB)_{i\ell}}

\newcommand{\bkl}{b_{k\ell}}
\newcommand{\btwoil}{b^{(2)}_{i\ell}}
\newcommand{\bll}{b_{\ell \ell}}

\newcommand{\matrixspace}{\;\;}
\newcommand{\ellone}{\ell_{1}}
\newcommand{\elltwo}{\ell_{2}}

\newcommand{\varphik}{\varphi_{k}}
\newcommand{\chik}{\chi_{k}}
\newcommand{\Phik}{\Phi_{k}}
\newcommand{\psik}{\psi_{k}}
\newcommand{\dr}{\beta^{n-k}}
\newcommand{\dn}{\beta^{-k}}
\newcommand{\betatwotilde}{\tilde{\rtwonk}}
\newcommand{\betathtilde}{\tilde{\ronenk}}
\newcommand{\betamink}{\beta^{-k}}

\newcommand{\rvarphik}{r(\varphi_{k})}    
\newcommand{\rrvarphik}{(r(\varphi_{k}))}    
\newcommand{\rchik}{r(\chi_{k})}    
\newcommand{\rrchik}{(r(\chi_{k}))}    
\newcommand{\svarphik}{s(\varphi_{k})}    
\newcommand{\ssvarphik}{(s(\varphi_{k}))}    
\newcommand{\schik}{s(\chi_{k})}    
\newcommand{\sschik}{(s(\chi_{k}))}    

\newcommand{\rpsik}{r(\psi_{k})}    
\newcommand{\rrpsik}{(r(\psi_{k}))} 
\newcommand{\spsik}{s(\psi_{k})}    
\newcommand{\sspsik}{(s(\psi_{k}))} 
\newcommand{\rpsinmk}{r(\psi_{n-k})}    
\newcommand{\rrpsinmk}{(r(\psi_{n-k}))} 
\newcommand{\spsinmk}{s(\psi_{n-k})}    
\newcommand{\sspsinmk}{(s(\psi_{n-k}))} 

\newcommand{\rvarphinmk}{r(\varphi_{n-k})}    
\newcommand{\rrvarphinmk}{(r(\varphi_{n-k}))}    
\newcommand{\rchinmk}{r(\chi_{n-k})}    
\newcommand{\rrchinmk}{(r(\chi_{n-k}))}    
\newcommand{\svarphinmk}{s(\varphi_{n-k})}    
\newcommand{\ssvarphinmk}{(s(\varphi_{n-k}))}    
\newcommand{\schinmk}{s(\chi_{n-k})}    
\newcommand{\sschinmk}{(s(\chi_{n-k}))}    
\newcommand{\stilde}{\tilde{s}}

\newcommand{\resh}{\rho}
\newcommand{\snk}{s}
\newcommand{\ronenk}{r_{1}}
\newcommand{\rtwonk}{r_{2}}
\newcommand{\dltaone}{\delta_{1}}
\newcommand{\dltatwo}{\delta_{2}}
\newcommand{\dltatld}{\tilde{\delta}}
\newcommand{\tauone}{\tau_{1}}
\newcommand{\tautwo}{\tau_{2}}
\newcommand{\xb}{\bar{x}}
\newcommand{\qone}{q_{1}}
\newcommand{\qtwo}{q_{2}}
\newcommand{\ub}{\bar{u}}
\newcommand{\cmm}{\complex^{m \times m}}
\newcommand{\opdsk}{\mathcal{O}}

\newcommand{\sqrtaone}{\sqrt{\a1}}
\newcommand{\sqrtatwo}{\sqrt{\atwo}}

\newcommand{\rhoone}{\rho_{1}}
\newcommand{\rhotwo}{\rho_{2}}
\newcommand{\xone}{x_{1}}
\newcommand{\xtwo}{x_{2}}
\newcommand{\xthr}{x_{3}}
\newcommand{\xfor}{x_{4}}
\newcommand{\xfiv}{x_{5}}

\newcommand{\lnorm}{\left \|}
\newcommand{\rnorm}{\right \|}
\newcommand{\lnrm}{\biggl | \biggl |}
\newcommand{\rnrm}{\biggr |\biggr |}

\newcommand{\recipp}{p^{\protect \#}}
\newcommand{\recipP}{P^{\protect \#}}

\newcommand{\varphione}{\varphi_{1}}
\newcommand{\varphitwo}{\varphi_{2}}

\newcommand{\psione}{\psi_{1}}
\newcommand{\psitwo}{\psi_{2}}

\newcommand{\absatwo}{|\atwo|}
\newcommand{\absqrtatwo}{\sqrt{|\atwo|}}
\newcommand{\absrhoone}{|\rhoone|}
\newcommand{\absrhotwo}{|\rhotwo|}
\newcommand{\sqrgamma}{\sqrt{\gamma}}
\newcommand{\aminb}{a-b}                 
\newcommand{\sqrtc}{\sqrt{c}}            
\newcommand{\sqrtabsc}{\sqrt{|c|}}            
\newcommand{\absalf}{|\alpha|}            
\newcommand{\rhoi}{\rho_{i}}                  
\newcommand{\sigmai}{\sigma_{i}}                  
\newcommand{\rhoj}{\rho_{j}}                  
\newcommand{\sigmaj}{\sigma_{j}}                  
\newcommand{\xij}{x_{ij}}                  
\newcommand{\yij}{y_{ij}}                  
\newcommand{\varphii}{\varphi_{i}}
\newcommand{\varphij}{\varphi_{j}}
\newcommand{\xstar}{x^{\ast}}
\newcommand{\uuone}{\mathcal{U}_{1}}
\newcommand{\uutwo}{\mathcal{U}_{2}}
\newcommand{\xkova}{\hat{x}}
\newcommand{\xkovalam}{\xkova_{max}(\ell-1)}

\begin{center}
\large
{\bf POLYNOMIAL EIGENVALUE BOUNDS FROM COMPANION FORMS}
\vskip 0.5cm
\normalsize
A. Melman \\
Department of Applied Mathematics \\
School of Engineering, Santa Clara University  \\
Santa Clara, CA 95053  \\
e-mail : amelman@scu.edu \\
\vskip 0.5cm
\end{center}

\begin{abstract}
We show how $\ell$-ifications, which are companion forms of matrix polynomials, namely, lower order matrix polynomials with the same eigenvalues 
as a given complex square matrix polynomial, can be used in combination with other recent results to produce eigenvalue bounds.
\vskip 0.15cm
{\bf AMS(MOS) subject classification :} 15A18, 47A46, 65F15
\end{abstract}

%
%
%
%

\section{Introduction}           
\label{introduction} 

Square matrix polynomials occur in polynomial eigenvalue problems, which consist of finding 
a nonzero complex eigenvector $v$ and a complex eigenvalue $z$ such that $P(z)v=0$, 
where $P$ is a matrix polynomial of the form
\bdis
\An  z^{n} + \Anmo z^{n-1} + \dots \Aone z + \Azero  \; ,
\edis
and $\Aj$ ($j=1,\dots,n$) are complex $m \times m$ matrices. If $A_{n}$ is singular then there are infinite eigenvalues, and if $A_{0}$ is singular then zero 
is an eigenvalue. 
There are $nm$ eigenvalues, including possibly infinite ones. The matrix polynomial is \emph{regular} if its determinant is not identically zero, and they are
the only ones we consider here.
The finite eigenvalues are the solutions 
of $\text{det} P(z)=0$. Polynomial eigenvalue problems can be found throughout many fields of engineering, a good overview of which can be found, e.g., 
in~\cite{BHMST}, \cite{TisseurMeerbergen}, and their references. The computation of polynomial eigenvalues is an active area of research.

As with the computation of zeros of a \emph{scalar} polynomial, whose coefficients are complex numbers, there are, broadly speaking, two approaches to compute
the eigenvalues of a matrix polynomial. One can either work with the matrix polynomial directly,
or one can compute the eigenvalues of its companion matrix. The disadvantage of the latter is that it can significantly increase the size of the matrix,
compared to the size of the coefficient matrices. 

Using the companion matrix amounts to a linearization of the matrix polynomial, i.e., the polynomial
eigenvalue problem has its order $n$ reduced and is replaced by a linear eigenvalue problem. 
In~\cite{DTDM} and~\cite{DTDVD} so-called $\ell$-ifications of a matrix polynomial are derived, which are, roughly speaking, matrix polynomials that are intermediate
between the given matrix polynomial and its linearization, i.e., they are spectrally equivalent matrix polynomials of a lower degree.
This also increases the size of the matrix coefficients, but to a lesser extent than would be the case when using the companion matrix. In analogy to the 
companion matrix, one might also call such $\ell$-ifications companion forms (\cite{DTDM}, \cite{DTDVD}), or companion matrix polynomials. 
The results in~\cite{DTDM} and~\cite{DTDVD} are very general and apply to rectangular matrix polynomials, whose elements can belong to any field. 

Our purpose here is to use the aforementioned $\ell$-ifications to construct bounds on the eigenvalues in the important special case of square complex 
matrix polynomials. To this end we will 
apply the generalization of a result by Cauchy from~\cite{HighamTisseur} (later also proved, in different ways, in~\cite{BiniNoferiniSharify} and~\cite{Melman_MatPol})
to different $\ell$-ifications (or companion forms) and then enhance those bounds with the generalization in~\cite{Melman_RS} of a result from~\cite{RS}.
As such, $\ell$-ification and enhancement of the bounds work hand in hand: the $\ell$-ification reduces the degree of the matrix polynomial, 
often significantly, and the enhancement prevents this from affecting the quality of the bounds. 

The paper is organized as follows. 
In Section~\ref{companion}, we define the $\ell$-ifications from~\cite{DTDM} and~\cite{DTDVD} for a square matrix polynomial, and derive the bounds
in Section~\ref{application}. 

%
%
%
%

\section{$\ell$-ifications of square complex matrix polynomials}
\label{companion}   

Throughout, we will denote the identity matrix by $I$ and the null matrix by $0$, without 
specifying their size, which should be clear from the context. A blank entry in a matrix denotes a zero or a zero block.
We define the exchange matrix $J$ as
\bdis
J = 
\begin{pmatrix}
   &          & 1  \\
   &  \iddots &    \\
1  &          &    \\
\end{pmatrix}
\; .
\edis
In addition, quantities with inadmissible indices are assumed to be zero, e.g.,
$b_{-1}=0$ if $b_{j}$ is not defined for negative integers.
We denote the reverse polynomial of $P$ by $P^{\protect \#}$, i.e., $P^{\protect \#}(z) = z^{n} P(1/z)$. 

Let us now consider a matrix polynomial $P(z)=\sum_{j=0}^{n} \Aj z^{j}$ with $A_{j} \in \complex^{m \times m}$.
If $n$ is divisible by a positive integer $k$, then, with $k < n$ and $q=\frac{n}{k}$, it was shown in~\cite[Section 4.4]{DTDVD}
that a strong $\ell$-ification of $P$, i.e., an $\ell$-ification that preserves the full eigenvalue structure of $P$, is given by  
$\sum_{j=0}^{q} J C_{j} J z^{j}$, where $J$ is the $km \times km$ exchange matrix, and the $k \times k$ block matrices 
$C_{j} \in \complex^{km \times km}$ are defined by
\begin{eqnarray}
& & C_{0} = 
\begin{pmatrix}
A_{(k-1)q} & A_{(k-2)q} & \dots  &  A_{q} & \Azero   \\
-I         &  0         &        &        &          \\
           &  -I        & \ddots  &        &          \\
           &            & \ddots &  0     &          \\
           &            &        &  -I    & 0        \\
\end{pmatrix}
\; ,
\;\;
C_{q} =
\begin{pmatrix}
\An   &       &      &        &           \\
      & I     &      &        &           \\
      &       & I    &        &           \\
      &       &      & \ddots &           \\
      &       &      &        & I         \\
\end{pmatrix}
\; ,
\nonumber \\
& & \nonumber \\
& & \text{and} \;\; C_{j} = 
\begin{pmatrix}
A_{j+(k-1)q} &  A_{j+(k-2)q}     &  \dots     & A_{j+q}   & A_{j}         \\
0            &  0                &  \dots     & 0         & 0             \\
\vdots       & \vdots            & \vdots     & \vdots    & \vdots        \\
0            &  0                & 0          & \dots     & 0             \\
\end{pmatrix}
\;\;\;\; (j=1,...,q-1)
\; .
\nonumber 
\end{eqnarray}
When $k=n$, then $q=1$,
\bdis
C_{0} = 
\begin{pmatrix}
\Anmo      & \Anmt      & \dots  & \Aone  & \Azero   \\
-I         &  0         &        &        &          \\
           &  -I        & \ddots  &        &          \\
           &            & \ddots &  0     &          \\
           &            &        &  -I    & 0        \\
\end{pmatrix}
\;\;\;\; \text{and} \;\;\;\;
C_{1} =
\begin{pmatrix}
\An   &       &      &        &           \\
      & I     &      &        &           \\
      &       & I    &        &           \\
      &       &      & \ddots &           \\
      &       &      &        & I         \\
\end{pmatrix}
\; .
\edis
Since $J^{2}=I$, the eigenvalues and eigenstructure of the $\ell$-ifications $\sum_{j=0}^{q} J C_{j} J z^{j}$ and $\sum_{j=0}^{q} C_{j} z^{j}$ are identical, 
and we prefer to work with the latter for reasons of convenience. The construction of such an $\ell$-ification in~\cite{DTDVD} is complicated because
of its generality and its strong result showing the equivalence of the eigenstructure of the $\ell$-ification and that of the given matrix polynomial. 

However, here we only need part of that result since we are only concerned with the eigenvalues themselves, not the eigenstructure. In addition, the matrix
polynomials we consider are square. 
This narrower focus allows us, in the theorem below, to provide a simpler path to $\ell$-ifications, avoiding the need 
of stepping through the relatively laborious process in~\cite{DTDVD}. 
It is based on a straightforward observation and on the following standard lemma, whose proof we defer to the appendix.
%
%
\begin{lemma}   
\label{Lemma_det}
For a positive integer $n \geq 3$, let $M_{j} \in \complex^{m \times m}$ for $j=1,...,n$ and $N_{j} \in \complex ^{m \times m}$ for $j=1,...,n-1$.
Then 
\bdis
\text{det} 
\begin{pmatrix}
\Mone & \Mtwo & \Mth   & \cdots & \Mn    \\
-I    & \None &        &        &        \\
      & -I    & \Ntwo  &        &        \\
      &       & \ddots & \ddots &        \\
      &       &        & -I     & \Nnmo  \\
\end{pmatrix}
= \text{det} \left ( \Mone \prod_{j=1}^{n-1} \Nj + \Mtwo \prod_{j=2}^{n-1} \Nj   + \dots + \Mnmo \Nnmo + \Mn \right ) \; ,
\edis
where the matrix multiplications are from the right with increasing index $j$.
\end{lemma}
The following theorem is a special case of the $\ell$-ification in~\cite[Section 4.4]{DTDVD}.   
%
%
\begin{theorem} 
\label{Theorem_companion}
Let $P(z)=\sum_{j=0}^{n} \Aj z^{j}$ and $Q(z)=\sum_{j=0}^{q} C_{j} z^{j}$, with the matrices $\Aj$ and $C_{j}$ defined as before.
If $n$ is divisible by a positive integer $k$, let $q=\frac{n}{k}$. Then 
the eigenvalues of $P$ and $Q$ coincide.
\end{theorem}
\prf
The proof is based on the observation that $P$ can be written as follows:
\begin{eqnarray}
\label{Pexpression} 
P(z) & = & \An z^{n} + \Anmo z^{n-1} + ... + \Azero   \nonumber \\
& = & z^{(k-1)\frac{n}{k}} \left ( \An z^{\frac{n}{k}} + \Anmo z^{\frac{n}{k} -1} + ... + A_{(k-1)\frac{n}{k}} \right ) \nonumber \\ 
& & \hskip 0.50cm + z^{(k-2)\frac{n}{k}} \left ( A_{(k-1)\frac{n}{k}-1} z^{\frac{n}{k}-1} + A_{(k-1)\frac{n}{k}-2} z^{\frac{n}{k}-2} 
               + ... + A_{(k-2)\frac{n}{k}} \right ) + \;\; ...     \nonumber \\
& & \hskip 0.50cm + z^{(k-j)\frac{n}{k}} \left ( A_{(k-(j-1))\frac{n}{k}-1} z^{\frac{n}{k}-1} + A_{(k-(j-1))\frac{n}{k}-2} z^{\frac{n}{k}-2} 
               + ... + A_{(k-j)\frac{n}{k}} \right ) + \;\; ...     \nonumber \\
& & \hskip 0.5cm + z^{\frac{n}{k}} \left ( A_{2\frac{n}{k}-1} z^{\frac{n}{k}-1} + A_{2\frac{n}{k}-2} z^{\frac{n}{k}-2} 
               + ... + A_{\frac{n}{k}} \right ) + A_{\frac{n}{k}-1} z^{\frac{n}{k}-1} + A_{\frac{n}{k}-2} z^{\frac{n}{k}-2} 
               + ... + \Azero  \; . \nonumber       \\
\end{eqnarray}
With $q=\frac{n}{k}$ and defining 
\beq
\label{M_def}
M_{i} (z)  = \sum_{j=0}^{q-1} A_{j+(k-i)q} z^{j} \;\; (i=1,...k) \; ,
\eeq   
the expression in~(\ref{Pexpression}) can be written as
\bdis
\hskip -0.25cm
P(z) = \sum_{j=0}^{n} \Aj z^{j} = \sum_{i=1}^{k} \lb \sum_{j=0}^{q-1} A_{j+(k-i)q} z^{j} \rb  z^{(k-i)q}
= \sum_{i=1}^{k} M_{i} (z)  z^{(k-i)q} \; .
\edis
Setting $M_{i}=M_{i}(z)$ for $i=1,...,k$ and $N_{i}=Iz^{q}$ for $i=1,...,k-1$ in Lemma~\ref{Lemma_det} yields        
\beq 
\label{detP_matrix}
\text{det}(P(z)) = \text{det} \lb \sum_{i=1}^{k} M_{i}(z) z^{(k-i)q} \rb 
=
\text{det} 
\begin{pmatrix}
\Mone & \Mtwo  & \Mth    & \cdots  & \Mk    \\
-I    & Iz^{q} &         &         &        \\
      & -I     & Iz^{q}  &         &        \\
      &        & \ddots  & \ddots  &        \\
      &        &         & -I      & Iz^{q} \\
\end{pmatrix}
\; .
\eeq
Using the definition of $M_{i}(z)$ in~(\ref{M_def}) to compare powers of $z$, it is now straightforward to verify from the right-hand side of~(\ref{detP_matrix})
that
\beq
\label{detPdetQ}
\text{det}(P(z)) = \text{det} \lb \sum_{j=0}^{q} C_{j} z^{j} \rb \; ,
\eeq  
where the matrices $C_{j}$ are as defined before.
This means that the finite eigenvalues of the two matrix polynomials $P$ and $Q$ coincide. 
On the other hand, we also have that 
$\text{det} \lb P^{\protect \#}(z) \rb = \text{det} \lb z^{n} P(1/z) \rb = z^{nm} \text{det} \lb P(1/z) \rb$, so that, with~(\ref{detPdetQ}), we obtain 
\bdis
\text{det} \lb P^{\protect \#}(z) \rb 
= z^{nm} \text{det} \lb P(1/z) \rb 
= z^{q(km)} \text{det} \lb Q(1/z) \rb 
= \text{det} \lb z^{q} Q(1/z) \rb 
= \text{det} \lb Q^{\protect \#}(z) \rb \; ,  
\edis
which implies that the infinite eigenvalues of $P$ and $Q$ also coincide. This completes the proof. 
\qed

We observe that, if $P$ is monic, then so is $Q$. Furthermore, $\text{det}(C_{0}) = \text{det}(\Azero)$ and $\text{det}(C_{q}) = \text{det}(\An)$, and, 
when $k=n$ and $\An=I$, the matrix polynomial $Q$
we just defined becomes $Iz + C_{0}$, so that $\text{det}(P)=0$ is equivalent to $\text{det} (Iz + C_{0})=0$, i.e., $-C_{0}$ is a companion matrix
of the matrix polynomial $P$ (see, e.g., \cite[Theorem 1.1]{GLR}), often called the Frobenbius companion matrix. 
Another special case, namely, for $k=2$, was derived in~\cite{Melman_MatPol}.

$\ell$-ifications can also be obtained from the companion matrix by permutation similarity transformations, but it would be a far more complicated procedure 
than the proof presented above. We therefore merely illustrate it here for the matrix polynomial 
$Iz^{4}+\Ath z^{3}+\Atwo z^{2}+\Aone z+ \Azero$, with $k=2$ and $q=2$.
A permutation similarity transformation of its companion matrix yields
\bdis
\begin{pmatrix}
I          &  0         &  0     &   0       \\
0          &  0         &  I     &   0       \\
0          &  I         &  0     &   0       \\
0          &  0         &  0     &   I       \\
\end{pmatrix}  
\begin{pmatrix}
-\Ath      & -\Atwo     & -\Aone & -\Azero   \\
I          &  0         &  0     &   0       \\
0          &  I         &  0     &   0       \\
0          &  0         &  I     &   0       \\
\end{pmatrix}  
\begin{pmatrix}
I          &  0         &  0     &   0       \\
0          &  0         &  I     &   0       \\
0          &  I         &  0     &   0       \\
0          &  0         &  0     &   I       \\
\end{pmatrix}  
=                   
\begin{pmatrix}
-\Ath      & -\Aone     & -\Atwo & -\Azero   \\
0          &  0         &  I     &   0       \\
I          &  0         &  0     &   0       \\
0          &  I         &  0     &   0       \\
\end{pmatrix}  
,
\edis
\normalsize
which is an $\ell$-ification of the matrix polynomial
\bdis
\begin{pmatrix}
I     &  0   \\
0     &  I   \\
\end{pmatrix}  
z^{2}
+
\begin{pmatrix}
\Ath      & \Aone  \\
0         &  0      \\
\end{pmatrix}  
z
+ 
\begin{pmatrix}
\Atwo & \Azero   \\
-I    &   0       \\
\end{pmatrix}  
\; ,
\edis
and this is precisely the companion form from Theorem~\ref{Theorem_companion} for this case. We remark that a formal proof along these
lines, beyond our scope here, would establish the equivalence of the eigenvalue structures as well.
\vskip 0.25cm

\noindent {\bf Example.} As an example consider
the matrix polynomial
\beq
\label{examplepol}
P(z) = A_{9}z^{9} + A_{8}z^{8} + A_{7}z^{7} + A_{6}z^{6} + A_{5}z^{5} + A_{4}z^{4} + A_{3}z^{3} + A_{2}z^{2} + A_{1}z + \Azero \; .
\eeq
Setting $k=3$, so that $q=3$, we obtain that the eigenvalues of $P$ are the same as those of 
\bdis
Q(z) =
  \begin{pmatrix} A_{9} & 0     & 0     \\ 0  & I & 0 \\ 0 & 0  & I \\ \end{pmatrix} z^{3} 
+ \begin{pmatrix} A_{8} & A_{5} & A_{2} \\ 0  & 0 & 0 \\ 0 & 0  & 0 \\ \end{pmatrix} z^{2}  
+ \begin{pmatrix} A_{7} & A_{4} & A_{1} \\ 0  & 0 & 0 \\ 0 & 0  & 0 \\ \end{pmatrix} z  
+ \begin{pmatrix} A_{6} & A_{3} & A_{0} \\ -I & 0 & 0 \\ 0 & -I & 0 \\ \end{pmatrix}   
\; . 
\edis

%
%
%
%

\section{Application: eigenvalue bounds} 
\label{application}
We now show that $\ell$-ification, in addition to its theoretical significance, also leads to bounds that do 
not appear to have been considered elsewhere. To this end, we combine the concept of $\ell$-ification with several other results.
We rely on the following generalization in~\cite{HighamTisseur} (later also in~\cite{BiniNoferiniSharify} and~\cite{Melman_MatPol}) of a classical
result of Cauchy, which states that all the eigenvalues of the square complex matrix polynomial $P(z) = \An z^{n} + A_{n-1}z^{n-1} + \dots + A_{1}z + A_{0}$, 
with $\An$ nonsingular, lie in $|z| \leq r$, where $r$ is the unique positive solution of
\beq
\label{Cauchyeq}
\|\An^{-1} \|^{-1} z^{n} - \|A_{n-1}\|z^{n-1} - \dots - \|A_{1}\|z - \|A_{0}\| = 0 \; , 
\eeq    
for any matrix norm. We call $r$ the \emph{Cauchy radius} of $P$. The polynomial equation~(\ref{Cauchyeq}) can easily be solved with any
standard root-finding method. We note that, for $\An=I$ and $n=1$, we obtain the well-known fact that $\|\Azero\|$
is an upper bound on the eigenvalues of $\Azero$. It is noteworthy that the Cauchy radius is typically among the best bounds for many problems,
as can be seen from~\cite{HighamTisseur}, where an exhaustive comparison of eigenvalue bounds was carried out. Consequently, we will use the Cauchy radius
as our benchmark for comparing bounds.

We also rely on an improvement of this Cauchy radius, which is the generalization in~\cite{Melman_RS} of a result for scalar polynomials in~\cite{RS}. 
Since this result seems less well-known and requires a few words of explanation, we restate it here. It is Theorem~2.2 in~\cite{Melman_RS}.
%
%
\begin{theorem}
\label{Theorem_genRS}
Let $P(z) = \sum_{j=0}^{n} A_{j} z^{j}$ be a regular matrix polynomial other than a matrix monomial, with $\Aj \in \complex^{m \times m}$ and 
$\An$ nonsingular. Denote by $i$ the smallest positive integer such that $A_{n-i}$ is not the null matrix, and 
define $T^{(L)}(z) = \lb \An z^{i} - \Anmi \rb P(z)$ and $T^{(R)}(z) = P(z) \lb \An z^{i} - \Anmi \rb$. 
If $\An\Anmi=\Anmi\An$ and $\|\An^{-2}\|^{-1} = \|\An\|\|\An^{-1}\|^{-1} $, then the Cauchy radii of $T^{(L)}$ and $T^{(R)}$ are not larger than the 
Cauchy radius of~$P$ when the same matrix norm is used for all radii.
\end{theorem}
Different norms yield different Cauchy radii, but it is hard to predict which norm delivers the smallest one.
The conditions $\|\An^{-2}\|^{-1} = \|\An\|\|\An^{-1}\|^{-1}$ and $\An\Anmi=\Anmi\An$ may appear to be restrictive, but they are always satisfied
when $\An=I$, which can be obtained by pre- or postmultiplication by $\An^{-1}$, and when $\|I\|=1$, which is the case for most standard norms (certainly for 
all subordinate norms).
The matrix $\An^{-1}$ needs to be computed in any case to apply the above generalized Cauchy result.
Lower bounds can be obtained analogously by applying the theorem to the reverse polynomial.

Theorem~\ref{Theorem_genRS} can be applied repeatedly to generate a nonincreasing sequence of Cauchy radii, with
the "left" and a "right" versions, in general, yielding different results.           
The improved Cauchy radii come at the cost of additional matrix multiplications and require the solution of real scalar polynomial equations 
of a degree higher than that of~$P$. It will therefore depend on the problem's properties (e.g., sparsity and/or symmetry) if this cost is justified. 
Already for moderately large matrix coefficients, the cost tends to be dominated by the matrix multiplications.

Clearly, by applying an $\ell$-ification to~$P$, one can lower the degree of the polynomial equation~(\ref{Cauchyeq}) that needs to be solved, which, unfortunately,
tends to decrease the quality of the resulting Cauchy radius. Fortunately, the improvement of the Cauchy radius from Theorem~\ref{Theorem_genRS}
mitigates this effect, thereby allowing $\ell$-ification to lower the degree of the polynomial without sacrificing the quality of the Cauchy radius.
Quite apart from this advantage, the bounds thus obtained do not seem to have appeared elsewhere in the literature.

Since we are considering upper bounds on the eigenvalues, we will assume from here on that our matrix polynomials have only finite eigenvalues, 
i.e., that their leading coefficients are nonsingular. Without loss of generality, we will therefore assume them to be monic. In addition, we will
only consider matrix norms for which $\|I\|=1$, such as the standard and often used $1$-norm, $\infty$-norm and $2$-norm.

Let us first revisit the above example before conducting a numerical comparison of Cauchy radii obtained through $\ell$-ification.
The Cauchy radius of $P$ in~(\ref{examplepol}) for a given matrix norm is the unique positve solution of 
\bdis
I z^{9} - \|A_{8}\|z^{8} - \|A_{7}\|z^{7} - \|A_{6}\|z^{6} - \|A_{5}\|z^{5} - \|A_{4}\|z^{4} - \|A_{3}\|z^{3} - \|A_{2}\|z^{2} - \|A_{1}\|z - \|\Azero\| = 0 \; .
\edis
As in the example, with $k=3$ and $q=3$, we obtain that the eigenvalues of $P$ are the same as those of the $\ell$-ification of $P$, given by
$Q(z) = Iz^{3}+C_{2}z^{2}+C_{1}z+C_{0}$, with
\bdis
C_{0} = \begin{pmatrix} A_{6} & A_{3} & A_{0} \\ -I & 0 & 0 \\ 0 & -I & 0 \\ \end{pmatrix}   
\;\; , \;\;
C_{1} = \begin{pmatrix} A_{7} & A_{4} & A_{1} \\ 0  & 0 & 0 \\ 0 & 0  & 0 \\ \end{pmatrix}   
\;\; , \;\; \text{and} \;\; 
C_{2}=\begin{pmatrix} A_{8} & A_{5} & A_{2} \\ 0  & 0 & 0 \\ 0 & 0  & 0 \\ \end{pmatrix} \; . 
\edis
Therefore, the Cauchy radius of $Q$, i.e., the unique positive solution of 
\bdis
z^{3} - \|C_{2}\|z^{2} - \|C_{1}\|z - \|C_{0}\| = 0 \; , 
\edis
is an upper bound on the eigenvalues of $P$. This Cauchy radius is obviously easier to compute than that of $P$ (it
can even be computed analytically), but it will, in general, be worse (i.e., larger) than that of $P$. However, it can be improved with Theorem~\ref{Theorem_genRS},
which states, when $C_{2} \neq 0$, that the Cauchy radii of $Q_{1}^{(L)}(z)=\lb Iz-C_{2} \rb Q(z)$ and $Q_{1}^{(R)}(z)=Q(z) \lb Iz-C_{2} \rb$ 
will be no worse than that of $Q$. In practice, they are often much better. 

It is instructive to compare $Q_{1}^{(L)}(z)=Iz^{4} + \lb C_{1} - C_{2}^{2} \rb z^{2} + \lb C_{0} - C_{2}C_{1} \rb z - C_{2} C_{0}$ and 
$Q_{1}^{(R)}(z)=Iz^{4} + \lb C_{1} - C_{2}^{2} \rb z^{2} + \lb C_{0} - C_{1}C_{2} \rb z - C_{0} C_{2}$. We observe that, in general, the degree of the polynomial
has been increased by one and that a zero coefficient has appeared immediately following the leading coefficient. A second application of the theorem
raises the degree to $6$ and adds an additional zero coefficient. More zero coefficients can appear, depending on the structure of the polynomial. 
The coefficient matrices of $Q_{1}^{(L)}$ are given by    
\begin{eqnarray*}
& & - C_{2} C_{0} = \begin{pmatrix} A_{5}-A_{8}A_{6} & A_{2}-A_{8}A_{3} & -A_{8} A_{0} \\ 0 & 0 & 0 \\ 0 & 0 & 0 \\ \end{pmatrix}   
\; , \\ 
& & C_{0} - C_{2}C_{1} = \begin{pmatrix} A_{6}-A_{8}A_{7} & A_{3}-A_{8}A_{4} & A_{0}-A_{8}A_{1} \\ -I  & 0 & 0 \\ 0 & -I  & 0 \\ \end{pmatrix}   
\; ,  \\
& & C_{1} - C_{2}^{2}=\begin{pmatrix} A_{7}-A_{8}^{2}& A_{4}-A_{8}A_{7} & A_{1}-A_{8}A_{2} \\ 0  & 0 & 0 \\ 0 & 0  & 0 \\ \end{pmatrix}  
\; ,  
\end{eqnarray*}
whereas those of $Q_{1}^{(R)}$ that are different are given by
\begin{eqnarray*}
& & - C_{0} C_{2} = \begin{pmatrix} -A_{6}A_{8} & -A_{6}A_{5} & -A_{6} A_{2} \\ A_{8} & A_{5} & A_{2} \\ 0 & 0 & 0 \\ \end{pmatrix}   
\; , \\ 
& & C_{0} - C_{1}C_{2} = \begin{pmatrix} A_{6}-A_{7}A_{8} & A_{3}-A_{7}A_{5} & A_{0}-A_{7}A_{2} \\ -I  & 0 & 0 \\ 0 & -I  & 0 \\ \end{pmatrix}   
\; .             
\end{eqnarray*}
Comparing these coefficient matrices, one notices that their norms could be quite different. For example, the constant coefficients $-C_{2}C_{0}$
and $-C_{0}C_{2}$ have a very different structure: if, e.g., it were the case that $A_{5} \approx A_{8}A_{6}$ and $A_{2} \approx A_{8}A_{3}$, then
$\|C_{2}C_{0}\|$ might be significantly smaller than $\|C_{0}C_{2}\|$ for most norms, leading to a smaller Cauchy radius. Similar situations
arise for other relations between the coefficients of $P$. Left and right multiplication, especially since it can be alternated in successive applications
of Theorem~\ref{Theorem_genRS}, therefore adds considerable flexibility.

If we set $k=1$ and $q=9$, we obtain a linearization of $P$, namely, $Iz+C_{P}$, where $-C_{P}$ is the companion matrix of $P$, and the Cauchy radius 
simply becomes $\|C_{P}\|$. In this case, Theorem~\ref{Theorem_genRS} states that the Cauchy radius of 
$(Iz-C_{P})(Iz+C_{P})=Iz^{2}-C_{P}^{2}$, namely, $\|C_{P}^{2}\|^{1/2}$, is not worse than that of $P$. This is not a surprise since 
the eigenvalues of $C_{P}^{2}$ are the squares of those of $C_{P}$, so that we obtain for any eigenvalue $\lambda$ of $C_{P}$:
\bdis
|\lambda^{2}| \leq \|C_{P}^{2}\| \leq \|C_{P}\|^{2} \Longrightarrow |\lambda| \leq \|C_{P}^{2}\|^{1/2} \leq \|C_{P} \| \; .
\edis
In fact, such reasoning leads to Gelfand's formula, wich states that $\rho(C) = \lim_{t \rightarrow +\infty} \|C^{t}\|^{1/t}$, where $\rho(C)$ 
is the spectral radius of the complex square matrix $C$. The matrix $C_{P}^{2}$ is given by 
\bdis
C_{P}^{2} = 
\begin{pmatrix} A_{8}^{2} - A_{7} &  \dots  & A_{8}A_{2}-A_{1} & A_{8}A_{1}-A_{0} & A_{8}A_{0}   \\ 
                          -A_{8}  &  \dots  & -A_{2}           &        A_{1}     &     -A_{0}   \\ 
                          I       &         &                  &                  &              \\ 
                                  &  \ddots &                  &                  &              \\ 
                                  &         &   I              &      0           &     0        \\ 
\end{pmatrix}
\; .
\edis
For a single application of the "left" version of Theorem~\ref{Theorem_genRS}, the number of operations necessary for the matrix multiplications 
and additions is very similar 
for both the linearization and the $3$-ification, although repeated application of Theorem~\ref{Theorem_genRS} favors the $3$-ification,
which, more or less, preserves the original sparsity of the coefficients. The latter is not the case for powers of $C_{P}$, which  quickly fill up the matrix.
The same is true for the computation of norms, which for $3$-ification requires fewer computations already for one application of the theorem. 
These observations continue to hold in general, as we will see below. 

The main setup cost of Theorem~\ref{Theorem_genRS} is determined by the matrix multiplications and the 
computation of the norms, although the latter's contribution is minor for the $1$-norm and the $\infty$-norm. 
In what follows, we therefore concentrate on estimating the former.
This cost is determined by the sparsity of the coefficients, which are initially 
sparse: only the top block-row is nonzero for the nonconstant coefficients and the constant coefficient has an additional identity
matrix in the lower left part. With each application of Theorem~\ref{Theorem_genRS}, they begin to fill up in a pattern that is difficult to predict,
which also makes it difficult to calculate the precise cost of the matrix multiplications. However, given a measure of sparsity in the form of the 
number of nonzero elements in the coefficient matrices, a rough but useful estimate can be derived, as was done in the appendix: 
when applying the theorem to an $\ell$-ification with given $k$ and $m$, and after applying Theorem~\ref{Theorem_genRS} once or more,
the computational effort to compute the matrix multiplications is proportional to 
\beq
\label{costestimate}
\dfrac{s^{2}}{\nu k m} \; ,
\eeq
where $s$ is the total number of nonzero elements in all coefficients of the corresponding $\ell$-ification,
while $\nu$ is the number of its nonzero coefficients, excluding the leading coefficient. 



Let us now consider a few numerical examples to illustrate the bounds derived from $\ell$-ification and to serve as a general guide for their use.
This usefulness depends on several factors: the sparsity and structure of the original matrix polynomials' coefficients,
and the computational cost one is willing to bear to obtain better bounds. 

As was previously mentioned, Cauchy radii are typically among the best bounds available for polynomial eigenvalues and outperform those based on singular
values or other explicit bounds (see~\cite{HighamTisseur}). 
They therefore provide a good benchmark. 

To compare the quality and computational cost of bounds obtained from repeatedly applying Theorem~\ref{Theorem_genRS} to different $\ell$-ifications of a 
matrix polynomial, we randomly generated three classes of $100$ monic matrix polynomials. For each case, two tables were produced: one table detailing 
the computational
effort involved in the matrix computations in units of the effort required to carry out one application with the "left" version 
of Theorem~\ref{Theorem_genRS} to the original
matrix polynomial, and the other table listing the average ratios of the Cauchy radii to the modulus of the largest eigenvalue, i.e., the closer these numbers 
are to one, the better they are. Between parentheses are the degrees of the real polynomial equations that need to be solved to obtain the Cauchy radius.  
Each successive row in the tables indicates an additional application of the "left" version of Theorem~\ref{Theorem_genRS}. There are three successive applications. 
The top row represents the Cauchy radii without applying the theorem; it requires no matrix multiplications. 
We have used the $1$-norm throughout, which delivers better Cauchy radii given the structure of the 
coefficient matrices. The definitions of the three classes and the accompanying comments follow.

\noindent \underline{\bf Class I}
\newline
This class consists of $4 \times 4$ matrix polynomials of degree $18$, where the elements of the non-leading coefficients have real and  
and imaginary parts that are uniformly randomly distributed on the interval $[-2,2]$. We then compared Cauchy radii
for $k=18,9,6,3,2,1$, corresponding to $q=1,2,3,6,9,18$, respectively.
The case $k=1$ gives the Cauchy radius of the given matrix polynomial itself, whereas $k=18$ corresponds to the companion matrix. 
The results can be seen in Table~\ref{table1} and Table~\ref{table2}.
We observe that, not surprisingly, higher values of $q$ deliver better bounds, but it  
is noteworthy that even low-order $\ell$-ifications can produce significantly improved bounds with just a few applications of Theorem~\ref{Theorem_genRS}.
For example, with roughly the same amount of work, a quartic polynomial (in the second row) 
delivers a result ($1.68$) that is not much worse than the result ($1.53$), which 
requires the solution of an equation of degree $19$ (in the same row).
The results for the companion matrix are better than for other $\ell$-ifications, but this is deceptive as it is clear from Table~\ref{table1} that this 
requires significantly more computation that sharply increases with each application of Theorem~\ref{Theorem_genRS}.
The computational cost was estimated using~(\ref{costestimate}), where the numbers of nonzero elements in the coefficient matrices were obtained numerically.
%
%
\begin{table}[H]
\begin{center}
\small           
\begin{tabular}{cccccc}
$q=1$         & $q=2$       & $q=3$      &  $q=6$    &  $q=9$      & $q=18$      \\
              &             &            &           &             &             \\
* (1)         &   * (2)     &  * (3)     & * (6)     &   * (9)     & * (18)      \\  
              &             &            &           &             &             \\
1.5 (2)       &   1.2 (3)   &  1.1 (4)   &  1.1 (7)  &   1.0 (10)  & 1.0 (19)    \\  
              &             &            &           &             &             \\
4.0 (4)       &   2.1 (5)   &  2.0 (6)   &  2.0 (9)  &   2.0 (12)  & 1.9 (21)    \\  
              &             &            &           &             &             \\
8.4 (8)       &   4.8 (8)   &  3.0 (9)   &  2.9 (12) &   2.9 (15)  & 2.9 (24)    \\  
\end{tabular}
\caption{Computational cost of repeated application of Theorem~\ref{Theorem_genRS}for $n=18$ and $m=4$.}                              
\label{table1}
\end{center}
\end{table}
\vskip -1.0cm 
%
%
\begin{table}[H]
\begin{center}
\small           
\begin{tabular}{cccccc}
$q=1$      &     $q=2$ &     $q=3$ &     $q=6$   & $q=9$      & $q=18$      \\
           &           &           &             &            &             \\
2.63 (1)   & 2.56 (2)  & 2.53 (3)  & 2.43 (6)    & 2.40 (9)   & 2.27 (18)   \\
           &           &           &             &            &             \\
1.80 (2)   & 1.71 (3)  & 1.68 (4)  & 1.63 (7)    & 1.59 (10)  & 1.53 (19)   \\
           &           &           &             &            &             \\
1.34 (4)   & 1.40 (5)  & 1.37 (6)  & 1.34 (9)    & 1.33 (12)  & 1.29 (21)   \\
           &           &           &             &            &             \\
1.15 (8)   & 1.34 (8)  & 1.31 (9)  & 1.29 (12)   & 1.27 (15)  & 1.24 (24)   \\
           &           &           &             &            &             \\
\end{tabular}
\vskip -0.45cm
\caption{Cauchy radius to maximum eigenvalue modulus ratios for $n=18$ and $m=4$.}                              
\label{table2}
\end{center}
\end{table}
\normalsize
\noindent \underline{\bf Class II}
\newline
This class consists of $100 \times 100$ matrix polynomials of degree $10$, where the elements of the non-leading coefficients have real and  
and imaginary parts that are uniformly randomly distributed on the interval $[-2,2]$. Here we compared Cauchy radii
for $k=10,5,2,1$, corresponding to $q=1,2,5,10$, respectively.
The results can be seen in Table~\ref{table3} and Table~\ref{table4}.
The Cauchy radii are much worse here than for the previous class of matrix polynomials, due to the larger size of the coefficient matrices.
The improvements are therefore more significant. It is somewhat remarkable that a polynomial of degree $5$ (in row 3) produces a bound (1.79)
that is very close to a bound (1.73), obtained with a polynomial of degree $16$ (in the bottom row) that, in fact, requires more work.
%
%
%
\begin{table}[H]
\begin{center}
\small           
\begin{tabular}{cccccc}
$q=1$         &   $q=2$     &  $q=5$      &  $q=10$       \\
              &             &             &               \\
* (1)         & * (2)       & * (5)       & * (10)        \\
              &             &             &               \\
1.0 (2)       & 1.0 (3)     & 1.0 (6)     & 1.0 (11)      \\
              &             &             &               \\
3.0 (4)       & 1.7 (5)     & 1.8 (8)     & 1.9 (13)      \\
              &             &             &               \\
7.0 (8)       & 4.2 (8)     & 2.7 (11)    & 2.8 (16)      \\
              &             &             &               \\
\end{tabular}
\vskip -0.45cm
\caption{Computational cost of repeated application of Theorem~\ref{Theorem_genRS} for $n=10$ and $m=100$.}                              
\label{table3}
\end{center}
\end{table}
\vskip -0.5cm 
%
%
\begin{table}[H]
\begin{center}
\small           
\begin{tabular}{cccc}
$q=1$      &   $q=2$   &  $q=5$    &  $q=10$      \\
           &           &           &              \\
9.89 (1)   & 9.84 (2)  & 9.72 (5)  & 9.66 (10)    \\
           &           &           &              \\
3.15 (2)   & 3.08 (3)  & 3.05 (6)  & 3.04 (11)    \\
           &           &           &              \\
1.77 (4)   & 1.79 (5)  & 1.78 (8)  & 1.77 (13)    \\
           &           &           &              \\
1.31 (8)   & 1.75 (8)  & 1.74 (11) &  1.73 (16)   \\
           &           &           &              \\
\end{tabular}
\vskip -0.45cm
\caption{Cauchy radius to maximum eigenvalue modulus ratios for $n=10$ and $m=100$.}                              
\label{table4}
\end{center}
\end{table}
\normalsize
\noindent \underline{\bf Class III}
\newline
Here we consider $10 \times 10$ matrix polynomials of degree $100$, where the elements of the non-leading coefficients have real and  
and imaginary parts that are uniformly randomly distributed on the interval $[-2,2]$. We compared Cauchy radii
for $k=100,50,25,20,10,5,2,1$, corresponding to $q=1,2,4,5,10,20,25,50,100$, respectively.
The results can be seen in Table~\ref{table5} and Table~\ref{table6}.
The quality of the Cauchy radii as bounds lie somewhere between the bounds in Class~I and Class~II since the same is true for the matrix sizes.
The degree of the polynomials is much higher than before, which makes it especially interesting that a quintic polynomial (in the third row)
delivers a bound (1.50) not much worse than one (1.39) obtained from a polynomial of degree $103$ (in the same row) for essentially the same amount of setup work.
%
%
\begin{table}[H]
\begin{center}
\small           
\begin{tabular}{ccccccccc}
$q=1$            &   $q=2$        &  $q=4$         &  $q=5$     &  $q=10$    &  $q=20$  &    $q=25$    &     $q=50$    &   $q=100$       \\
                 &                &                &            &            &          &              &               &                 \\
* (1)            & * (2)          & * (4)          & * (5)      &   * (10)   &  * (20)  &   * (25)     &   * (50)      &   * (100)       \\
                 &                &                &            &            &          &              &               &                 \\
1.2 (2)          & 1.1 (3)        & 1.0 (5)        & 1.1 (6)    & 1.0 (11)   & 1.0 (21) &  1.0 (26)    &   1.0 (51)    & 1.0 (101)       \\
                 &                &                &            &            &          &              &               &                 \\
3.4 (4)          & 1.8 (5)        & 1.9 (7)        & 1.9 (8)    & 1.9 (13)   & 2.0 (23) &  2.0 (28)    &   2.0 (53)    & 2.0 (103)       \\
                 &                &                &            &            &          &              &               &                 \\
7.6 (8)          & 4.4 (8)        & 2.8 (10)       & 2.9 (11)   & 2.9 (16)   & 2.9 (26) &  2.9 (31)    &   3.0 (56)    & 3.0 (106)       \\
                 &                &                &            &            &          &              &               &                 \\
\end{tabular}
\caption{Computational cost of repeated application of Theorem~\ref{Theorem_genRS} for $n=100$ and $m=10$.}                              
\label{table5}
\end{center}
\end{table}
\normalsize
%
%
\begin{table}[H]
\begin{center}
\small           
\begin{tabular}{ccccccccc}
$q=1$            &   $q=2$        &  $q=4$         &  $q=5$     &  $q=10$    &  $q=20$   &    $q=25$    &     $q=50$    &   $q=100$       \\
                 &                &                &            &            &           &              &               &                 \\
3.75 (1)         & 3.71 (2)       & 3.65 (4)       & 3.61 (5)   & 3.54 (10)  & 3.48 (20) & 3.46 (25)    & 3.37 (50)     & 3.27 (100)      \\ 
                 &                &                &            &            &           &              &               &                 \\
2.09 (2)         & 2.01 (3)       & 1.98 (5)       & 1.98 (6)   & 1.94 (11)  & 1.90 (21) & 1.89 (26)    & 1.85 (51)     & 1.80 (101)      \\
                 &                &                &            &            &           &              &               &                 \\
1.45 (4)         & 1.50 (5)       & 1.48 (7)       & 1.48 (8)   & 1.46 (13)  & 1.44 (23) & 1.44 (28)    & 1.41 (53)     & 1.39 (103)       \\
                 &                &                &            &            &           &              &               &                 \\
1.19 (8)         & 1.44 (8)       & 1.43 (10)      & 1.42 (11)  & 1.40 (16)  & 1.39 (26) & 1.38 (31)    & 1.36 (56)     & 1.34 (106)       \\
                 &                &                &            &            &           &              &               &                 \\
\end{tabular}
\caption{Cauchy radius to maximum eigenvalue modulus ratios for $n=100$ and $m=10$.}                              
\label{table6}
\end{center}
\end{table}
\normalsize

\vskip 0.5cm
\noindent {\bf Conclusion}
\newline
We have combined the recent concept of $\ell$-ification with other results into a framework for the construction of bounds on the moduli of polynomial 
eigenvalues. The computational cost and the quality of the resulting bounds were investigated at the hand of randomly generated classes of matrix polynomials 
of varying degrees and sizes, providing information that can be helpful in the choice of bounds.
Within this framework, there is room for many variations that cannot all be addressed here. For example, "left" and "right" versions of Theorem~\ref{Theorem_genRS}
can be applied or different matrix norms can be used. In addition, similarity transformations could further enhance the bounds   
by lowering the norms of the coefficient matrices, and polynomials could be multiplied by an appropriate power of $z$ to be able to apply a certain 
$\ell$-ification. The usefulness of such variations will depend on the properties of the matrix polynomials concerned.

%
%
%
%

\section{Appendix} 
\label{appendix}             

\noindent \underline{\bf Proof of Lemma~\ref{Lemma_det}}
\newline
We recall that Lemma~\ref{Lemma_det} states that, for $n \geq 3$ and $M_{j} , N_{j} \in \complex^{m \times m}$,
\bdis
\text{det} 
\begin{pmatrix}
\Mone & \Mtwo & \Mth   & \cdots & \Mn    \\
-I    & \None &        &        &        \\
      & -I    & \Ntwo  &        &        \\
      &       & \ddots & \ddots &        \\
      &       &        & -I     & \Nnmo  \\
\end{pmatrix}
= \text{det} \left ( \Mone \prod_{j=1}^{n-1} \Nj + \Mtwo \prod_{j=2}^{n-1} \Nj   + \dots + \Mnmo \Nnmo + \Mn \right ) \; ,
\edis
\normalsize
where the matrix multiplications are from the right with increasing index $j$.

Before we prove the lemma, we observe that (see, e.g., \cite[p.27]{HJ}), for complex matrices $A,B,C,D$ of the same size
with $CD=DC$, 
\bdis
\text{det} 
\begin{pmatrix}
A     & B        \\
C     & D        \\
\end{pmatrix}
= \text{det} \lb AD - BC \rb \; .
\edis
\prf
The proof is by induction. We first show that the lemma holds for a $3 \times 3$ block matrix, for which we have 
\begin{eqnarray}
\text{det} 
\begin{pmatrix}
\Mone & \Mtwo & \Mth  \\
-I    & \None & 0     \\
0     & -I    & \Ntwo \\
\end{pmatrix}
& = &
\text{det} 
\left (
\begin{pmatrix}
\Mone & \Mtwo & \Mth  \\
-I    & \None & 0     \\
0     & -I    & \Ntwo \\
\end{pmatrix}
\begin{pmatrix}
I     & 0     & 0      \\
0     & I     & \Ntwo  \\
0     & 0     & I      \\
\end{pmatrix}
\right )
\nonumber  \\
& = & 
\text{det} 
\begin{pmatrix}
\Mone & \Mtwo & \Mtwo \Ntwo + \Mth  \\
-I    & \None & \None \Ntwo         \\
0     &  -I   & 0                   \\
\end{pmatrix}
\label{basis1}   
\\
& = & 
\text{det} 
\begin{pmatrix}
\Mone & \Mtwo \Ntwo + \Mth  \\
-I    & \None \Ntwo         \\
\end{pmatrix}
\label{basis2}
\\
& = &
\text{det} \lb \Mone \None \Ntwo + \Mtwo \Ntwo + \Mth \rb \; , 
\label{basis3}
\end{eqnarray}
\normalsize
where~(\ref{basis2}) is obtained using Laplace expansion on the $-I$ block in the bottom block row in~(\ref{basis1}), and~(\ref{basis3}) follows from the fact that 
$\None \Ntwo$ and $-I$ commute. 
This establishes the induction basis. For the induction step, assume that the lemma holds for a $n \times n$ block matrix. For a $(n+1) \times (n+1)$ block matrix
we then have
\bdis
\text{det} 
\begin{pmatrix}
\Mone & \Mtwo & \Mth  & \cdots & \Mn & \Mnpo \\
-I    & \None &       &        &           &       \\
      & -I    & \Ntwo &        &           &       \\
      &       & -I    & \ddots &           &       \\
      &       &       & \ddots & \Nnmo     &       \\
      &       &       &        & -I        & \Nn   \\
\end{pmatrix}
\edis
\begin{eqnarray*}
& = &
\text{det} 
\left (
\begin{pmatrix}
\Mone & \Mtwo & \Mth  & \cdots & \Mn & \Mnpo \\
-I    & \None &       &        &           &       \\
      & -I    & \Ntwo &        &           &       \\
      &       & -I    & \ddots &           &       \\
      &       &       & \ddots & \Nnmo     &       \\
      &       &       &        & -I        & \Nn   \\
\end{pmatrix}
\begin{pmatrix}
I     &       &      &        &           &       \\
      & I     &      &        &           &       \\
      &       & I    &        &           &       \\
      &       &      & \ddots &           &       \\
      &       &      &        & I         & \Nn   \\
      &       &      &        &           &  I    \\
\end{pmatrix}
\right )
\\
& & \\
& = & 
\text{det} 
\begin{pmatrix}
\Mone & \Mtwo & \Mth  & \cdots & \Mn       & \Mn \Nn + \Mnpo \\
-I    & \None &       &        &           &                 \\
      & -I    & \Ntwo &        &           &                 \\
      &       & -I    & \ddots &           &                 \\
      &       &       & \ddots & \Nnmo     & \Nnmo \Nn       \\
      &       &       &        & -I        & 0               \\
\end{pmatrix}
\; ,
\end{eqnarray*}
\normalsize
which, by using Laplace expansion on the $-I$ block in the bottom block row, is equal to
\bdis
\text{det}
\begin{pmatrix}
\Mone & \Mtwo & \Mth  & \cdots & \Mnmo     & \Mn \Nn + \Mnpo  \\
-I    & \None &       &        &           &                  \\
      & -I    & \Ntwo &        &           &                  \\
      &       & -I    & \ddots &           &                  \\
      &       &       & \ddots & \Nnmt     &                  \\
      &       &       &        & -I        & \Nnmo \Nn        \\
\end{pmatrix}
\; .
\edis
\normalsize
Finally, using the induction hypothesis, this becomes
\bdis
\text{det} \left ( \Mone \lb \prod_{j=1}^{n-2} \Nj \rb \Nnmo \Nn + \Mtwo \lb \prod_{j=2}^{n-2} \Nj \rb \Nnmo \Nn  
+ \dots + \Mnmo \Nnmo \Nn +  \Mn \Nn + \Mnpo \right ) ,
\edis 
which completes the proof. \qed 

\noindent \underline{\bf Matrix multiplication cost estimate}
\newline
As we saw before, each application of Theorem~\ref{Theorem_genRS} adds at least one zero coefficient immediately after the leading 
coefficient of the matrix polynomial to which is applied, while increasing its degree. 
For a particular value of $k$, we denote by $s$ the total number of nonzero elements in all coefficients of the corresponding $\ell$-ification,
and by $\nu$ the number of nonzero coefficients, excluding the leading coefficient.
We will assume that these nonzero elements are more or less equally distributed among the coefficients, i.e., each coefficient contains 
on average $s/\nu$ nonzero elements. That means that each coefficient contains, on average, $\lb s/\nu \rb / \lb km^{2} \rb$ block-rows of  
$m \times m$ matrices. To apply Theorem~\ref{Theorem_genRS}, one such coefficient multiplies all the others, which requires
\bdis
k \, \lb \dfrac{s/\nu}{km^{2}}  \rb^{2}
\edis
$m \times m$ matrix multiplications per nonzero nonleading coefficient. Since each of those these matrix multiplications requires a number of 
operations proportional 
to $m^{3}$, the total computational cost is proportional to
\bdis
k \, \lb \dfrac{s/\nu}{km^{2}}  \rb^{2} m^{3} \nu = \dfrac{s^{2}}{\nu k m} \; . 
\edis  
Although crude, this expression provides an adequate way to compare the computational cost of the different bounds obtained from different
values of $k$ and repeated applications of Theorem~\ref{Theorem_genRS}.

\end{document}